\ifx\shlhetal\undefinedcontrolsequenc\let\shlhetal\relax\fi

\input amstex
\NoBlackBoxes
\documentstyle {amsppt}
\topmatter
\title {Finite Canonization \\
Sh564} \endtitle
\author {Saharon Shelah \thanks{\null\newline
I thank Alice Leonhardt for the beautiful typing \null\newline
Written 4/June/94 - Publ.No.564 \null\newline
Latest Revision - Aug/4/95} \endthanks} \endauthor
\affil {Institute of Mathematics \\
The Hebrew University \\
Jerusalem, Israel
\medskip
Rutgers University \\
Department of Mathematics \\
New Brunswick, NJ USA} \endaffil

\abstract    The canonization theorem says that for given $m,n$ for 
some $m^*$ (the first one is called $ER(n;m)$) we have 
\medskip
\roster
\item "{{}}"  for every function $f$ with domain $[{1,\dotsc,m^*}]^n$,
for some $A \in [{1,\dotsc,m^*}]^m$, the question of when the equality
$f({i_1,\dotsc,i_n}) = f({j_1,\dotsc,j_n})$ (where \newline
$i_1 < \cdots < i_n$ and
$j_1 < \cdots j_n$ are from $A$) holds has the simplest answer: for 
some $v \subseteq \{1,\dotsc,n\}$ the equality holds iff $\dsize 
\bigwedge_{\ell \in v} i_\ell = j_\ell$.
\endroster
\medskip

We improve the bound on $ER(n,m)$ so that fixing $n$ the number of
exponentiation needed to calculate $ER(n,m)$ is best possible. 
\endabstract
\endtopmatter
\document

\expandafter\ifx\csname bib4plain.tex\endcsname\relax
  \expandafter\gdef\csname bib4plain.tex\endcsname{}
\else \message{Hey!  Apparently you were trying to \string twice.   This does not make sense.}
\errmessage{Please edit your file (probably \jobname.tex) and remove
any duplicate ``\string\input'' lines} \fi

\def\renewcommand{\newcommand}	       
\edef\cite{\the\catcode`@}%
\catcode`@ = 11
\let\@oldatcatcode = \cite
\chardef\@letter = 11
\chardef\@other = 12
%
%
%
%
\def\@innerdef#1#2{\edef#1{\expandafter\noexpand\csname #2\endcsname}}%
%
%
\@innerdef\@innernewcount{newcount}%
\@innerdef\@innernewdimen{newdimen}%
\@innerdef\@innernewif{newif}%
\@innerdef\@innernewwrite{newwrite}%
%
%
%
\def\@gobble#1{}%
%
%
%
\ifx\inputlineno\@undefined
   \let\@linenumber = \empty 
\else
   \def\@linenumber{\the\inputlineno:\space}%
\fi
%
%
%
\def\@futurenonspacelet#1{\def\cs{#1}%
   \afterassignment\@stepone\let\@nexttoken=
}%
\begingroup 
\def\\{\global\let\@stoken= }%
\\ 
\endgroup
\def\@stepone{\expandafter\futurelet\cs\@steptwo}%
\def\@steptwo{\expandafter\ifx\cs\@stoken\let\@@next=\@stepthree
   \else\let\@@next=\@nexttoken\fi \@@next}%
\def\@stepthree{\afterassignment\@stepone\let\@@next= }%
%
%
%
\def\@getoptionalarg#1{%
   \let\@optionaltemp = #1%
   \let\@optionalnext = \relax
   \@futurenonspacelet\@optionalnext\@bracketcheck
}%
%
%
\def\@bracketcheck{%
   \ifx [\@optionalnext
      \expandafter\@@getoptionalarg
   \else
      \let\@optionalarg = \empty
      \expandafter\@optionaltemp
   \fi
}%
\def\@@getoptionalarg[#1]{%
   \def\@optionalarg{#1}%
   \@optionaltemp
}%
%
%
%
\def\@nnil{\@nil}%
\def\@fornoop#1\@@#2#3{}%
\def\@for#1:=#2\do#3{%
   \edef\@fortmp{#2}%
   \ifx\@fortmp\empty \else
      \expandafter\@forloop#2,\@nil,\@nil\@@#1{#3}%
   \fi
}%
\def\@forloop#1,#2,#3\@@#4#5{\def#4{#1}\ifx #4\@nnil \else
       #5\def#4{#2}\ifx #4\@nnil \else#5\@iforloop #3\@@#4{#5}\fi\fi
}%
\def\@iforloop#1,#2\@@#3#4{\def#3{#1}\ifx #3\@nnil
       \let\@nextwhile=\@fornoop \else
      #4\relax\let\@nextwhile=\@iforloop\fi\@nextwhile#2\@@#3{#4}%
}%
%
%
%
\@innernewif\if@fileexists
\def\@testfileexistence{\@getoptionalarg\@finishtestfileexistence}%
\def\@finishtestfileexistence#1{%
   \begingroup
      \def\extension{#1}%
      \immediate\openin0 =
         \ifx\@optionalarg\empty\jobname\else\@optionalarg\fi
         \ifx\extension\empty \else .#1\fi
         \space
      \ifeof 0
         \global\@fileexistsfalse
      \else
         \global\@fileexiststrue
      \fi
      \immediate\closein0
   \endgroup
}%
%
%
%
%
\def\bibliographystyle#1{%
   \@readauxfile
   \@writeaux{\string\bibstyle{#1}}%
}%
\let\bibstyle = \@gobble
%
%
\let\bblfilebasename = \jobname
\def\bibliography#1{%
   \@readauxfile
   \@writeaux{\string\bibdata{#1}}%
   \@testfileexistence[\bblfilebasename]{bbl}%
   \if@fileexists
      \nobreak
      \@readbblfile
   \fi
}%
\let\bibdata = \@gobble
%
%
\def\nocite#1{%
   \@readauxfile
   \@writeaux{\string\citation{#1}}%
}%
\@innernewif\if@notfirstcitation
%
%
\def\cite{\@getoptionalarg\@cite}%
%
%
\def\@cite#1{%
   \let\@citenotetext = \@optionalarg
   \printcitestart
   \nocite{#1}%
   \@notfirstcitationfalse
   \@for \@citation :=#1\do
   {%
      \expandafter\@onecitation\@citation\@@
   }%
   \ifx\empty\@citenotetext\else
      \printcitenote{\@citenotetext}%
   \fi
   \printcitefinish
}%
\def\@onecitation#1\@@{%
   \if@notfirstcitation
      \printbetweencitations
   \fi
   \expandafter \ifx \csname\@citelabel{#1}\endcsname \relax
      \if@citewarning
         \message{\@linenumber Undefined citation `#1'.}%
      \fi
      \expandafter\gdef\csname\@citelabel{#1}\endcsname{%
\strut
\vadjust{\vskip-\dp\strutbox
\vbox to 0pt{\vss\parindent0cm \leftskip=\hsize 
\advance\leftskip3mm
\advance\hsize 4cm\strut\openup-4pt 
\rightskip 0cm plus 1cm minus 0.5cm ?  #1 ?\strut}}
         {\tt
            \escapechar = -1
            \nobreak\hskip0pt
            \expandafter\string\csname#1\endcsname
            \nobreak\hskip0pt
         }%
      }%
   \fi
   \csname\@citelabel{#1}\endcsname
   \@notfirstcitationtrue
}%
%
%
\def\@citelabel#1{b@#1}%
%
%
\def\@citedef#1#2{\expandafter\gdef\csname\@citelabel{#1}\endcsname{#2}}%
%
%
%
\def\@readbblfile{%
   \ifx\@itemnum\@undefined
      \@innernewcount\@itemnum
   \fi
   \begingroup
      \def\begin##1##2{%
         \setbox0 = \hbox{\biblabelcontents{##2}}%
         \biblabelwidth = \wd0
      }%
      \def\end##1{}
      %
      %
      \@itemnum = 0
      \def\bibitem{\@getoptionalarg\@bibitem}%
      \def\@bibitem{%
         \ifx\@optionalarg\empty
            \expandafter\@numberedbibitem
         \else
            \expandafter\@alphabibitem
         \fi
      }%
      \def\@alphabibitem##1{%
         \expandafter \xdef\csname\@citelabel{##1}\endcsname {\@optionalarg}%
         \ifx\biblabelprecontents\@undefined
            \let\biblabelprecontents = \relax
         \fi
         \ifx\biblabelpostcontents\@undefined
            \let\biblabelpostcontents = \hss
         \fi
         \@finishbibitem{##1}%
      }%
      \def\@numberedbibitem##1{%
         \advance\@itemnum by 1
         \expandafter \xdef\csname\@citelabel{##1}\endcsname{\number\@itemnum}%
         \ifx\biblabelprecontents\@undefined
            \let\biblabelprecontents = \hss
         \fi
         \ifx\biblabelpostcontents\@undefined
            \let\biblabelpostcontents = \relax
         \fi
         \@finishbibitem{##1}%
      }%
      \def\@finishbibitem##1{%
         \biblabelprint{\csname\@citelabel{##1}\endcsname}%
         \@writeaux{\string\@citedef{##1}{\csname\@citelabel{##1}\endcsname}}%
         \ignorespaces
      }%
      %
      %
      \let\em = \bblem
      \let\newblock = \bblnewblock
      \let\sc = \bblsc
      \frenchspacing
      \clubpenalty = 4000 \widowpenalty = 4000
      \tolerance = 10000 \hfuzz = .5pt
      \everypar = {\hangindent = \biblabelwidth
                      \advance\hangindent by \biblabelextraspace}%
      \bblrm
      \parskip = 1.5ex plus .5ex minus .5ex
      \biblabelextraspace = .5em
      \bblhook
      \input \bblfilebasename.bbl
   \endgroup
}%
%
%
\@innernewdimen\biblabelwidth
\@innernewdimen\biblabelextraspace
%
%
%
\def\biblabelprint#1{%
   \noindent
   \hbox to \biblabelwidth{%
      \biblabelprecontents
      \biblabelcontents{#1}%
      \biblabelpostcontents
   }%
   \kern\biblabelextraspace
}%
%
%
%
\def\biblabelcontents#1{{\bblrm [#1]}}%
%
%
\def\bblrm{\rm}%
%
%
\def\bblem{\it}%
%
%
\def\bblsc{\ifx\@scfont\@undefined
              \font\@scfont = cmcsc10
           \fi
           \@scfont
}%
%
%
\def\bblnewblock{\hskip .11em plus .33em minus .07em }%
%
%
\let\bblhook = \empty
%
%
%
\def\printcitestart{[}
\def\printcitefinish{]}
\def\printbetweencitations{, }
\def\printcitenote#1{, #1}
%
%
%
\let\citation = \@gobble
%
%
%
\@innernewcount\@numparams
%
%
\def\newcommand#1{%
   \def\@commandname{#1}%
   \@getoptionalarg\@continuenewcommand
}%
%
%
\def\@continuenewcommand{%
   \@numparams = \ifx\@optionalarg\empty 0\else\@optionalarg \fi \relax
   \@newcommand
}%
%
%
\def\@newcommand#1{%
   \def\@startdef{\expandafter\edef\@commandname}%
   \ifnum\@numparams=0
      \let\@paramdef = \empty
   \else
      \ifnum\@numparams>9
         \errmessage{\the\@numparams\space is too many parameters}%
      \else
         \ifnum\@numparams<0
            \errmessage{\the\@numparams\space is too few parameters}%
         \else
            \edef\@paramdef{%
               \ifcase\@numparams
                  \empty  No arguments.
               \or ####1%
               \or ####1####2%
               \or ####1####2####3%
               \or ####1####2####3####4%
               \or ####1####2####3####4####5%
               \or ####1####2####3####4####5####6%
               \or ####1####2####3####4####5####6####7%
               \or ####1####2####3####4####5####6####7####8%
               \or ####1####2####3####4####5####6####7####8####9%
               \fi
            }%
         \fi
      \fi
   \fi
   \expandafter\@startdef\@paramdef{#1}%
}%
%
%
%
%
\def\@readauxfile{%
   \if@auxfiledone \else 
      \global\@auxfiledonetrue
      \@testfileexistence{aux}%
      \if@fileexists
         \begingroup
            \endlinechar = -1
            \catcode`@ = 11
            \input \jobname.aux
         \endgroup
      \else
         \message{\@undefinedmessage}%
         \global\@citewarningfalse
      \fi
      \immediate\openout\@auxfile = \jobname.aux
   \fi
}%
%
%
\newif\if@auxfiledone
\ifx\noauxfile\@undefined \else \@auxfiledonetrue\fi
%
%
%
%
\@innernewwrite\@auxfile
\def\@writeaux#1{\ifx\noauxfile\@undefined \write\@auxfile{#1}\fi}%
%
%
%
\ifx\@undefinedmessage\@undefined
   \def\@undefinedmessage{No .aux file; I won't give you warnings about
                          undefined citations.}%
\fi
%
%
\@innernewif\if@citewarning
\ifx\noauxfile\@undefined \@citewarningtrue\fi
%
%
%
\catcode`@ = \@oldatcatcode

\newpage

\head {\S0 Introduction} \endhead
\bigskip

On Ramsey theory see the book Graham Rotchild Spencer \cite{GrRoSp}.
This paper is self-contained.

The canonical Ramsey theorem was originally proved by Erdos and Rado,
so the relevant number is called $ER(n,m)$.  See \cite{ErRa},
\cite{Ra86} and more in the work of Galvin.  The theorem states that
if $m$ and $n$ are given, and $f$ is an $n$-place function on a set
$A$ of size $\ge ER(n,m)$, then there is an $A' \in [A]^m$ such that
$f$ is canonical on $A'$.  That is, for some $v \subseteq \{
1,\dotsc,n\}$ and for every $i_1 < \cdots < i_n \in A'$ and $j_1 <
\cdots < j_n \in A$

$$
f(i_1,\dotsc,i_n) = f(j_1,\dotsc,j_n) \Leftrightarrow \dsize \bigwedge
_{\ell \in v} i_\ell = j_\ell.
$$
\medskip

\noindent

Galvin got in the early seventies by the probability method a lower
bound which appeared in \cite{ErSp},p.30, $ER(2;m) \ge (m + o(1))^m$.
\newline
Lefmann and Rodl \cite{LeRo93} proved

$$
2^{cm^2} < ER(2;m) \le 2^{2^{c_1^{m^3}}}.
$$

\medskip
\noindent
Lefmann and Rodl \cite{LeRo94} proved:
\medskip
\roster
\item "{(i)}"  $2^{c_2m^2} \le ER(2;m) \le 2^{c^*_2(m^2 \text{ log }m)}$
\item "{(ii)}"  $\beth_n(c_km^2) \le ER(n+1;m) \le \beth_{n+1}
\left ( c^*_k {\frac {m^{2k-1}}{\text{log }m}} \right)$.
\endroster
\medskip

\noindent
See more on this in \cite{LeRo94} and below for the definition of
$\beth_n$. \newline
We thank Nesetril for telling us the problem; which for us was finding the
right number of exponents (i.e. the subscript for $\beth$ in $(ii)$ above)
in $ER(n;m)$ (for a fixed $n$).  We prove here that
this number is $n$. \newline
Why is the number of exponentiations best possible?  Let $r^n_t(m)$ be the
first $r$ such that: $r \rightarrow (m)^n_t$, now trivially
$ER(n;m) \ge r^n_t(m)$ when $m$ is not too small, and $r^n_t(m)$ needs
$n-1$ exponentiations when $t$ is not too small.
\newpage

\head {\S1 The finitary canonization lemma} \endhead
\bigskip

\demo{Notation}  $\Bbb R,\Bbb N$ are the set of reals and natural numbers
respectively.  The letters $k,\ell,m,n$ will be used to denote natural
numbers, as well as $i,j,\alpha, \beta,\gamma,\zeta,\xi$.  We let
$\varepsilon$ be a real (usually positive). 
\newline
If $A$ is a set,
$$
[A]^n = \{ u \subseteq A:|u| = n \}.
$$
\medskip

\noindent
We call finite subsets $u,v$ of $\Bbb N$ neighbors if:
$$
|u| = |v|,|u \backslash v| = 1
$$
and
$$
[k \in u \backslash v,\ell \in v \backslash u,m \in u \cap v \Rightarrow
k < m \equiv \ell < m].
$$
\medskip

\noindent
For $m \in \Bbb N$, we let $[m] = \{ 1,\dotsc,m\}$.

\noindent
For a set $A$ of natural numbers and $i \in \Bbb N$, $A<i$  means
$(\forall j \in A)(j < i)$. We similarly define $i < A$.

With $i,A$ as above
$$
A_{> i} \text{ denotes the set } \{j \in A:j > i\}.
$$
\medskip

\noindent
We use the convention that $A_{> \text{sup }\emptyset}$ is $A$.

\noindent
Let $\beth_n(m)$ be defined by induction on $n:\beth_0(m) = m$ and
$\beth_{n+1}(m) = 2^{\beth_n(m)}$. \newline
Usually, $c_i$ denotes a constant.
\enddemo
\bigskip

\proclaim{1.1 Lemma(Finitary Canonization)}  Assume $n$ is given, then there
is a constant $c$ computable from $n$, such that if $m$ is large enough:
\newline
If $f$ is an $n$-place function from $[m^\otimes] 
= \{1,\dotsc,m^\otimes\}$ and $m^\otimes > \beth_{n-1}(cm^{8(2n-1)})$
\newline
\underbar{then} for some $A' \in [\{ 1,\dotsc,m^\otimes\}]^m,f$ is canonical
on $A'$; i.e. for some $v \subseteq \{1,\dotsc,n\}$ for every $i_1 < \cdots
< i_n$ from $A'$ and $j_1 < \cdots < j_n$ from $A'$, we have
$$
f(i_1,\dotsc,i_n) = f(j_1,\dotsc,j_n) \Leftrightarrow \dsize \bigwedge
_{\ell \in v} i_\ell = j_\ell.
$$
\endproclaim
\bigskip

\noindent
The proof is broken into several claims.
\medskip

\noindent
\underbar{Explanation} of our proof. \newline

We inductively on $n^* = n^\otimes,\dotsc,1$ decrease the set to
$A_{n^*}$ while increasing the amount of ``partial homogeneity", i.e. 
conditions close to: results of computing $f$ on an $n$-tuple from
$A_{n^*}$ are not dependent on the last $p = n^\otimes - n^*$ 
members of the $n$-tuple.
Having gone down from $n^\otimes$ to $n^*$, we want that: if
$u_1,u_2 \in [A_{n^*}]^n$ are neighbors differing in the $\ell$-th place
element only then: if $\ell < n^*$, the truth value of $f(u_1) = f(u_2)$ 
depends on the first $n^*$ elements of $u_1$ and $u_2$ only; if $\ell > n^*$
the truth value of $f(u_1) = f(u_2)$ depends on the first $n^*$ elements of
$u_1$ only.  Lastly if $\ell = n^*$, it is little more complicated to control
this; but the truth value is monotonic and we introduce certain functions, 
(the $h$'s) which express this.  Arriving
to $n^*=1$ we eliminate the $h$'s (decreasing a little) so we get the 
sufficiency
of the condition for equality, but we still have the necessity only for 
$u_1,u_2$ which are neighbors.  Then by
random choice (as in \cite{Sh:37}), we get the necessity for all 
pairs of sets.  The earlier
steps cost essentially one exponentiation each, the last two cost only taking
a power.
\bigskip

\proclaim{1.2 Claim}  Assume
\medskip
\roster
\item "{$(*)_0$}"  $m \ge 2^{(1+\varepsilon)c_1(m^*)^{n^*}}$
\newline
$t > 0,n^* > 1,k(*) > 0 \, (c_1$ is defined in the proof from $k(*),n^*$)
\newline
and $m^*$ is large enough (relative to $1/\varepsilon,t,k(*),n^*)$
\item "{$(*)_1$}"  $A \subseteq \Bbb N,|A| > m$ \newline
$f_k$ a function with domain $[A]^{n^*}$ for $k < k(*)$, \newline
$h_k$ is a function from $[A]^{n^*}$ to $\Bbb N$ for $k < k(*)$, and \newline
$g$ is a function with domain $[A]^{n^*}$ such that $\text{Rang}(g)$ has
cardinality $\le t$.
\endroster
\medskip

\noindent
\underbar{Then} we can find $A^*,j^*$ such that:
\medskip
\roster
\item "{$(*)_2$}"  $A^* \subseteq A,|A^*| > m^*$ and $j^* \in
A_{> \text{ sup}(A^*)}$ and we have: \newline
if $k < k(*),u \in [A^*]^{n^*-1}$ and $v \in [A^*]^{n^*-1}$, \underbar{then}
\medskip
{\roster
\itemitem{ $(\alpha)$ }  if $u,v$ are neighbors, then 
for all $i \in A^*_{>\text{ sup}(u \cup \nu)}$ we have
$$
f_k(u \cup \{ i \}) = f_k(v \cup \{ i \}) \Leftrightarrow 
f_k(u \cup \{ j^*\}) = f_k(v \cup \{ j^*\})
$$
\medskip
\itemitem{ $(\beta)$ }   if $u=v$ then for every 
$i_0 < i_1$ from $A^*_{> \text{ sup}(u)}$ we have 
\footnote{This is used later to define the 
$h_k$ for the ``next step".} 
$$
f_k(u \cup \{ i_0 \}) = f_k(u \cup \{ i_1\}) \Leftrightarrow
f_k(u \cup \{ i_0 \}) = f_k(u \cup \{ j^* \})
$$
\medskip

\itemitem{ $(\gamma)$ }  for all $i \in A^*_{>\text{ sup}(u)}$
$$
g(u \cup \{ i \}) = g(u \cup \{ j^* \})
$$
\medskip

\itemitem { $(\delta)$ }  either for all $i \in A^*_{>\text{ sup}(u)} \cup
\{ j^* \}$ we have
$$
h_k(u) \ge i
$$
\noindent
\underbar{or} for all $i \in A^*_{>\text{ sup}(u)} \cup \{ j^* \}$ we have
$$
h_k(u) < i.
$$
\endroster}
\endroster
\endproclaim
\bigskip

\remark{1.2A Remark}  1) We could have also related $f_{k_1}(u),f_{k_2}(u)$
for various $k_1,k_2$, this would not have influenced the bounds.
\endremark
\bigskip

\demo{Proof}  Standard ramification. For $B \subseteq A$ we define an
equivalence relation $E_B$ on $A_{>\text{ sup}(B)}$ as follows. We let: 
\newline 
$i_0 \, E_B \, i_1$ \underbar{iff} $i_0,i_1 \in A_{>\text{ sup}(B)}$
and for every $u,v \in [B]^{n^*-1},d \in \text{Rang}(g)$, $w \in
[B]^{n^*}$ and $k < k(*)$ the truth value of the following is the same
for $\ell \in \{ 0,1 \}$: 
\medskip
\roster
\item "{$(\alpha)$}"  $f_k(u \cup \{ i_\ell \}) = f_k(v \cup \{ i_\ell \})$
if $u,v$ are neighbors
\item "{$(\beta)$}"  $f_k(u \cup \{ i_\ell \}) = f_k(w)$ if $u = w \backslash
\{\text{max}(w)\}$
\item "{$(\gamma)$}"  $g(u \cup \{ i_\ell \}) = d$
\item "{$(\delta)$}"  $h_k(u) \ge i_\ell$.
\endroster
\medskip

\noindent
Clearly $E_B$ is an equivalence relation and $E_\emptyset$ is the equality (as
$n^* > 1)$. \newline
For $i \in A_{> \text{ sup}(B)}$ we let $i/E_B$ denote the 
equivalence class of $i$ via $E_B$. \newline
Note that if $B \subseteq B^*$, then $i/E_{B^*} \subseteq i/E_B$.

\noindent
We now define a tree $T$ by defining by induction on $\ell \in \Bbb N$ 
objects $t_{\le \ell},\le_\ell$ and \newline
$\langle A_i:i \in t_{\le \ell} \rangle$ such that:
\medskip
\roster
\item "{(a)}"  $(t_{\le_\ell},\le_\ell)$ is a tree,$t_{\le_\ell}$ a subset
of $A,\le_\ell$ a partial order on $t_\ell$ such that for every
$x \in t_{\le_\ell},\{y:y \le_\ell x\}$ is linearly ordered
\item "{(b)}"  $t_{\le_\ell} \subseteq t_{\le_{\ell + 1}}$ and 
$\le_{\ell + 1} \restriction t_{\le_\ell} = \le_\ell$
\item "{(c)}"  $t_{\le 0} = \{ \text{min}(A)\},A_{\text{min}(A)} = 
A_{> \text{ min}(A)}$
\item "{(d)}"  $t_{\le(\ell + 1)} \backslash t_{\le \ell}$ is the
$(\ell + 1)$-th level of $(t_{\le (\ell + 1)},\le_{\ell + 1})$
\item "{(e)}"  if $i_0 <_\ell \, i_1 <_\ell \cdots <_\ell \, i_\ell \in 
t_{\le_\ell}$ (so $\{ i_0,\dotsc,i_\ell\}$ is a branch) then
{\roster
\itemitem { $(\alpha)$ }  $A_{i_\ell} = i_\ell/E_{\{i_0,\dotsc,i_{\ell-1}\}}$
\itemitem { $(\beta)$ }   the set of immediate successors of $i_\ell$ in
$(t_{\le(\ell + 1)},\le_{\ell + 1})$ is \newline
$Y_{i_\ell} =: \{ \text{min}(j/E_{\{i_0,i_1,\dotsc,i_\ell\}}):j \in
A_{i_\ell} \text{ but } j \ne i_\ell\}$.
\endroster}
\endroster
\medskip

\noindent
This is straight.  Let $t_\ell = t_{\le \ell} \backslash \dsize \bigcup
_{m < \ell} t_{\le m}$ and $T = \dsize \bigcup_\ell t_{\le \ell}$. \newline
Note also that $i \le_\ell j \Rightarrow i \le j$ and that
\medskip
\roster
\item "{$\bigotimes$}"  if we consider the definition of 
$E_{\{i:i \le_\ell j \}}$  restricted just to $A_j \backslash \{ j\}$ 
we may restrict ourselves:
for clause $(\alpha)$ only to the $u,v \in [\{ i:i \le_\ell j\}]
^{n^*-1}$ with \newline
$\text{max}(u \cup v) = j$, and for clause $(\beta)$ only to
those \newline
$u \in [\{ i:i \le_\ell j\}]^{n^*-1},w \in [\{ i:i \le_\ell j\}]^{n^*}$ with
$\text{max}(w) = j$. \newline
For $(\gamma)$ and $(\delta)$ we may assume $\text{max}(u) = i_\ell$.
\endroster
\medskip

\noindent
Now it is easy to see that
\medskip
\roster
\item "{$(*)_3$}"  $A = \dsize \bigcup_\ell t_\ell$
\item "{$(*)_4$}"  if $j \in t_\ell$ then the number of immediate
successors of $j$ in $(t_{\le \ell + 1},\le_{\ell + 1})$ (necessarily they
are all in $t_{\ell + 1}$) is at most

$$
\left( 2^{{\binom \ell{n^*-1}}(n^*-1)} \right)^{k(*)} \times
\left( 2^{\binom \ell{n^*-1}} \right)^{k(*)} \times t^{\binom \ell{n^*-2}}
\times \left( \tbinom \ell{n^*-2} \cdot k(*) + 1 \right).
$$
\endroster
\medskip

\noindent
[Why this inequality?  The four terms in the product correspond to the four
clauses $(\alpha), (\beta), (\gamma), (\delta)$ in the definition of
$E_B$ for the branch $B = \{ i_0,\dotsc,i_\ell = j\}$ of \newline
$(t_{\le \ell},\le)$.  The power $k(*)$ in the first two terms comes 
from dealing
with $f_k$ for each $k < k(*)$ and ``$2$ to the power $x$" as we have $x$
choices of yes/no.  Now the first term comes from counting the possible
$u \cup v$ (from clause $(\alpha)$).  At the first glance their number is
$|[\{i_0,\dotsc,i_\ell\}]^{n^*}|$ as being neighbors each 
with $n^*-1$ elements they
have together $n^*$ elements, but by $\bigotimes$ we can 
restrict ourselves to the
case $i_\ell \in u \cup v$, so we have to consider $|[\{i_0,\dotsc,
i_{\ell - 1}\}]^{n^*-1}| = \binom \ell{n^*-1}$ sets $u \cup v$; then we 
have to
choose $u \cup \nu \backslash (u \cap v)$ (as we do not need to distinguish
between $(u,v)$ and $(v,u)$).  As $u,v$ are neighbors we have $n^*-1$ possible
choices (as the two members of $(u \cup v) \backslash (u \cap v)$ 
are successive members of $u \cup v$ under the natural order).

For the second term, we should consider $u,w$ as in clause $(\beta)$, and so
as \newline
$u = w \backslash \{\text{max}(w)\}$ we know $w$  gives all the
information, and by $\bigotimes$ above \newline
$\text{max}(w) = i_\ell$, so the
number of possibilities is $\binom \ell{n^*-1}$. \newline

For the third term we have a choice of one from $t(=|\text{Rang}(g)|)$ for
each \newline
$u \in [\{i_0,\dotsc,i_\ell\}]^{n^*-1}$, but again by $\bigotimes$, with
$\text{max}(u) = i_\ell$, so the number is $\binom \ell{n^*-2}$.  \newline

Lastly, in the fourth term the number of questions ``$h_k(u) \ge i$" is again
\newline
$\binom \ell{n^*-2} \cdot k(*)$, but by the properties of linear orders 
there are $\binom \ell{n^*-2} \cdot k(*) + 1$ possible answers.  So 
$(*)_4$ really holds]. \newline
Clearly (with $c_0 = k(*)/(n^*-2)! + k(*)/(n^*-1)! + {\frac{\text{log}
(\ell^{n^* -1} \cdot k(*))}{\ell^{n^*-1}}}$ and 
$c^0_1 = c_0/n^*$)
\medskip
\roster
\item "{$(*)_5$}"  $\left( 2^{\binom \ell{n^*-1} \times n^*-1} \right)^{k(*)}
\times \left(2^{\binom \ell{n^*-1}} \right)^{k(*)} \times
t^{\binom \ell{n^*-2}} \times \left( \tbinom \ell{n^*-2} \cdot k(*) + 1 
\right)$ \newline
$\le 2^{k(*)\ell^{n^*-1}/(n^*-2)!} \times 2^{k(*) \cdot \ell^{n^*-1}/
(n^*-1)!} \times 2^{\text{log}(t)\ell^{n^*-2}/(n^*-2)!}$ \newline
$\times \ell^{n^*-1} \cdot k(*) \le 2^{c_0\ell^{n^*-1}(1+\varepsilon)}$.
\endroster
\medskip

\noindent
(any positive $\varepsilon$, for $\ell$ large enough). \newline
So (for some constant $c^2_0$)
\medskip
\roster
\item "{$(*)_6$}"  $|t_{\ell+1}| \le c^2_0 \dsize \prod^\ell_{p=1}
2^{(1 + \varepsilon)c_0 p^{n^*-1}} = c^2_0\cdot
2^{(1 + \varepsilon)c_0 \dsize \sum^\ell_{p=1} p^{n^*-1}}$ \newline
$\le c^2_0\cdot 2^{(1 + \varepsilon)c_0(\ell + 1)^{n^*}/n^*} = 
c^2_0\cdot 2^{(1 + \varepsilon) c^0_1(\ell+1)^{n^*}}$.
\endroster
\medskip

\noindent
But
\medskip
\roster
\item "{$\bigoplus$}"  if $a_p \ge 0,a_p \le a_{p+1}$ and $p \ge \ell^* 
\Rightarrow 2a_p \le
a_{p+1}$ then \newline
$\dsize \sum^\ell_{p=0} a_p \le 2a_\ell + \dsize \sum_{p \le \ell^*} a_p$ 
\endroster
\medskip

\noindent
hence (possibly increasing $\varepsilon$, which means for $(*)_5$ using
large $\ell$)
\roster
\item "{$(*)_7$}"   $|t_{\le(\ell + 1)}| \le c^3_0 +
\dsize \sum^{\ell + 1}_{p=0} c^2_0\cdot 2^{(1 + \varepsilon)c^0_1(p +
1)^{n^*}} \le c^4_0\cdot 2^{(1 + \varepsilon)c^0_1(\ell + 1)^{n^*} +
1} \le$ \newline  
$2^{(1 + \varepsilon)c_1(\ell + 1)^{n^*}}$.
\endroster
\medskip 

\noindent
So (increasing $\varepsilon$ slightly)
\medskip
\roster
\item "{$(*)_8$}"   $|t_{\le m^*}| \le m < |A|$
\endroster
\medskip

\noindent
so there is a $j^* \in t_{m^*+1}$.  Let
$A^* = \{ i:i <_{m^*+1}j^*\}$ (so $|A^*| = m^*+1)$, then \newline
$A^*,j^*$ are as required. \hfill$\square_{1.2}$
\enddemo
\bigskip

\proclaim{1.3 Claim}  Assume
\medskip
\roster
\item "{$(*)_9(a)$}"  $n^\otimes \ge n^* \ge 1,k(*) > 0$
\item "{$(b)$}"  we have the function $m(-)$ satisfying $m(n+1) \ge 
2^{(1+ \varepsilon)c_1m(n)^{n+1}}$ for \newline
$n \in [n^*,n^\otimes)$
\item "{$(c)$}"  $t > 0$ and $m(n^*)$ is large enough relative to $k(*),
n^*,c_1,1/\varepsilon$.
\endroster
\medskip
\roster
\item "{$(*)_{10}$}"  $A \subseteq \Bbb N,|A| \ge m(n^\otimes + 1)$, $g$ is a 
function
with domain $[A]^{n^\otimes}$ and $\text{range}$ with $\le t$ members;
$f_k$ is a function with domain $[A]^{n^\otimes}$ (for $k < k(*)$), 
and for simplicity
${\Cal P}(\{ 0,1,\dotsc,n^\otimes - 1\}) \cap \text{ Rang}(f_k) = \emptyset$.
\endroster
\medskip

\noindent
\underbar{Then} we can find $A' \in [A]^{m(n^*)+1}$ and $j^*_\ell \in A$
for $\ell \in [n^*,n^\otimes)$ satisfying \newline
$A' < j^*_{n^*} < j^*_{n^*+1} < \dots$, and functions
$g',g_k,h_k(k < k(*))$ with domain $[A']^{n^*}$ such that (letting
$w^* = \{ j^*_\ell:\ell \in [n^*,n^\otimes)\}$):
\medskip
\roster
\item "{$(*)_{11}$}"  for all $u \in [A']^{n^*}$
{\roster
\itemitem{ (a) } for $w \in [A'_{> \text{ sup}(u)}]^{n^\otimes-n^*}$ we
have \newline
$g(u \cup w) = g'(u) = g(u \cup w^*)$
\itemitem{ (b) }  for $k < k(*)$ we have $h_k(u) \in \Bbb N$ and
$g_k(u) \in \{ v:v \subseteq (n^*,n^\otimes) \}$
\itemitem{ (c) }  \underbar{if} $w_1,w_2 \in [A'_{> \text{ sup}(u)}]
^{n^\otimes -n^*}$ and $k < k(*)$ and: (note:$|u| = n^*$) \newline
$\{i \in w_1:|u| + |i \cap w_1| \in g_k(u)\} = \{ i \in w_2:|u| + |i \cap
w_2| \in g_k(u)\}$ \newline
and $[\text{min}(w_1 \cup w_2) < h_k(u) \Rightarrow \text{ min}(w_1) =
\text{ min}(w_2)]$ \underbar{then} \newline
$f_k(u \cup w_1) = f_k(u \cup w_2)$.
\itemitem{ (d) }  Assume $k < k(*),w_1 \cup w_2 \cup \{i,j \} \subseteq A',
u < w_1 < i < j < w_2$ and \newline
$|w_1 \cup w_2| = n^\otimes - n^* - 1$: \newline
\medskip

$\qquad (i)$  if $w_1 \ne \emptyset$ \underbar{then}
$$
f_k(u \cup w_1 \cup \{ i \} \cup w_2) = f_k(u \cup w_1 \cup \{ j \} \cup
w_2) \Leftrightarrow |u \cup w_1| \notin g_k(u)
$$
\medskip

$\qquad (ii)$ if $w_1 = \emptyset$ \underbar{then}
$$
f_k(u \cup \{ i \} \cup w_2) = f_k(u \cup \{ j \} \cup w_2) 
\Leftrightarrow h_k(u) \le i.
$$
\itemitem{ (e) }  for $k < k(*)$ and neighbors 
$u_0,u_1 \in [A']^{n^*}$ and \newline
$w \in [A'_{> \text{ max}(u_0 \cup u_1)}]^{n^\otimes - n^*}$ we have:
$$
f_k(u_0 \cup w) = f_k(u_1 \cup w) \text{ iff }
f_k(u_0 \cup w^*) = f_k(u_1 \cup w^*).
$$
\endroster}
\endroster
\endproclaim
\bigskip

\remark{Remark}  1) Note particularly clause ($d$).  
So $g_k(u)$ is intended to
be like the $v$ in 1.1, only fixing an initial segment of both
$\{ i_\ell:\ell < n^\otimes\}$ and $\{ j_\ell:\ell < n^\otimes\}$ as $u$.
But whereas the equality demand in clause $(d)$ is as expected, the
non-equality demand is weaker: only for neighbors. \newline 
2) Note that we can in some clauses above replace $A'$ by $A' \cup W^*$.
\endremark
\bigskip

\demo{Proof}  We prove this by induction on $n^\otimes - n^*$.  If it is zero,
the conclusion is trivial.
\newline
\noindent
Use the induction hypothesis with $n^\otimes,n^*+1,f_k,(k < k(*)),g$
now standing for $n^\otimes,n^*,f_k,(k < k(*)),g$ in the induction hypothesis.
We get $A' \in [A]^{m(n^*+1)+1}$ and
functions $g',g_k,h_k$ (for $k < k(*)$) and $j^*_\ell$ for $\ell \in
[n^* + 1,n^\otimes)$ satisfying $(*)_{11}$ of Claim 1.3.
Now we apply 1.2 to $m = m(n^* + 1), A', g^\otimes, f^\otimes_k,
h^\otimes_k(k < k(*))$ where
we define the function $g^\otimes$ with domain $[A']^{n^*+1}$ by
$g^\otimes(u) = \langle g'(u),g_k(u):k < k(*) \rangle,h^\otimes_k = h_k$
and the function
$f^\otimes_k$ with domain $[A']^{n^*+1}$ is defined by

$$
f^\otimes_k(u) = f_k(u \cup \{j^*_\ell:\ell \in [n^*+1,n^\otimes)\}).
$$
\medskip

\noindent
We get there $A^* \in [A']^{m(n^*) + 1}$ and $j^* \in (A')_{> \text{ sup }
A^*}$.  Let $j_{n^*} =: j^*$.  Now we have to define $h_k$ with domain
$[A^*]^{n^*}$ (for $k < k(*)$).  For $u \in [A^*]^{n^*}$ let

$$
B^k_u =: \{ i \in A^*_{> \text{ sup}(u)}:f^\otimes_k(u \cup \{ i \}) \ne
f^\otimes_k(u \cup \{ j^* \} \}.
$$
\medskip

\noindent
By clause $(\beta)$ of Claim 1.2, $B^k_i$ is an initial segment of
$A^*_{> \text{ sup}(u)}$.  Let $h_k(u) = \text{ max}(B^k_u) + 1$.

Lastly, for $u \in [A^*]^{n^*}$ we have to define $g_k(u)$.  By 1.2
$(\delta)$, the answer to \newline
``$h^\otimes_k(u \cup \{j\}) < j_{n^*}$" does not
depend on $u$ and on $j \in A^*_{>\text{ sup } u}$.  Let $g^\otimes_k$ be the
``old" $g_k$ (with domain $[A']^{n^*+1})$ and let

$$
g_k(u) = \cases g^\otimes_k(u \cup \{j_{n^*}\}) \quad &\text{ if } \quad
h^\otimes_k(u \cup \{j\}) < j_{n^*} \\
g^\otimes_k(n \cup \{ j_{n^*}\}) \cup \{ j_{n^*}\} \quad &\text{ otherwise}.
\endcases
$$
\medskip

\noindent
Now $A^*,g_k,h_k,j^*_n,j^*_{n^*+1},\dots$ are as required.
\hfill$\square_{1.3}$
\enddemo
\bigskip

\proclaim{1.4 Claim}  1) Assume $m(1) \ge (k(*) \cdot m(0))^{k(*)+1}$ and
$A' \subseteq \Bbb N,|A'| \ge m(1)$ and for
$k < k(*),h_k$ is a function from $A'$ into $\Bbb N,h_k(i) \ge i$.

Then we can find $A'' \subseteq A',|A''| \ge m(0)$ such that
\medskip
\roster
\item "{$(*)_{12}$}"  for each $k < k(*)$ we have:

$$
\text{either } \qquad (\forall i,j \in A'')[i < j \Rightarrow h_k(i) \ge j]
$$
$$
\text{or } \qquad \,\,(\forall i,j \in A'')[i < j \Rightarrow h_k(i) < j].
$$
\endroster
\medskip

\noindent
2)  If $m(1) > d^km(0)^{2^{k(*)}},A \subseteq \Bbb N,|A| > m(1),g_k$ is a
function from $A$ to $\{ 1,\dotsc,d\}$, and $f_k$ is a function from $A$ to
$\Bbb N$ for $k < k(*)$ then we can find $A' \subseteq A,|A'| > m(0)$ such
that:
\medskip
\roster
\item "{$\bigotimes$}"  for each $k,f_k \restriction A'$ is constant or
one to one and $g_k \restriction A'$ is constant.
\endroster
\endproclaim
\bigskip

\demo{Proof}  1) We can find $A_1 \subseteq A',|A_1| > m(1)/k(*)^{k(*)}$ such
that for all $i,j \in A_1$,

$$
\ell,k < k(*) \Rightarrow [h_\ell(i) \le h_k(i) \equiv h_\ell(j) \le
h_k(j)].
$$

\noindent
So without loss of generality 
\medskip
\roster
\item "{$(*)$}"  $\ell < k < k(*) \and i \in A_1 \Rightarrow h_\ell(i) \le
h_k(i)$.
\endroster
\medskip

\noindent
By renaming we can assume $A_1 = \{1,2,\dotsc,m(0)^{k(*)+1}\}$. \newline
Now if for some $\ell,0 < \ell \le m(0)^{k(*)+1} - m(0)$, and \newline
$(\forall \alpha)\biggl(\alpha \in [\ell,\ell + m(0)) \Rightarrow 
h_0(\alpha) \ge \ell + m(0) \biggr)$ then $A'' = [\ell,\ell + m(0))$ is as 
required for all $h_\ell$ by $(*)$. \newline
If not, then we can find $\alpha_\ell \in [1,m(0)^{k(*)+1})$ for
$\ell = 1,\dotsc,m(0)^{k(*)}$, strictly increasing with $\ell$ such that 
$h_0(\alpha_\ell) < \alpha_{\ell + 1}$.  We repeat the argument for 
$h_1$, etc. \newline
2) Also easy.  \hfill$\square_{1.4}$
\enddemo
\bigskip

\remark{Remark}  We can use $m(1) > k(*)! \cdot m(0)^{k(*)+1}$ instead.
The only point is the choice of $A$.
\endremark
\bigskip

\proclaim{1.5 Claim}  Assume we have the assumptions of 1.3.
If we first apply 1.3 getting $A'$ and then aply 1.4 to get 
$A'' \subseteq A'$ such that for each
$k$ and $u \in [A'']^{n^*}$ either \newline
$h_k(u) \le \text{ min}\{ \ell \in A'':u < \ell\}$ or 
$h_k(u) > \text{ max}(A'')$ (we assume now $n^*=1$ so $u = \{j\}$), and in
addition

\medskip
\roster
\item "{$(*)_{13}$}" $m^{2n^\otimes-n^*} \cdot (2n^\otimes - n^*)
{\binom {2(n^\otimes-n^*)-1}{n^\otimes-n^*}} \cdot k(*) \le |A''|$
\endroster
\medskip

\noindent
\underbar{then} there is $A^* \in [A'']^m$ such that (in addition to 
$(*)_{11}(a) - (e) + (*)_{12}$) we have
\medskip
\roster
\item "{$(*)_{14}$}"  for all $u \in [A^*]^{n^*}$ 
{\roster
\itemitem{ $(f)$ }  if $w_1,w_2 \in [A^*_{> \text{ sup}(u)}]^{n^\otimes-n^*},
k < k(*)$ \underbar{then}
\endroster}
\endroster
$$
\align
f_k(u \cup w_1) = &f_k(u \cup w_2) \Leftrightarrow
 \{ i \in w_1:|i \cap w_1| + |u| \in g_k(u) \} = \\  
  &\{ i \in w_2:|i \cap w_2| + |u| \in g_k(u)\}.
\endalign
$$
\endproclaim
\bigskip

\remark{Remark}  Here we are rectifying the gap between the equality 
($(*)_{11}(d)$) and the inequality ($(*)_{11}(e)$) demand.
\endremark
\bigskip

\demo{Proof}  First note that
\medskip
\roster
\item "{$(*)_{15}$}"  for all $u \in [A'']^{n^*}$ the implication 
$\Leftarrow$ holds. \newline
[why? just use clause ($c$) of $(*)_{11}$ of Claim 1.3].
\endroster
\medskip

\noindent
So we are left with proving $\Rightarrow$.
\enddemo
\bigskip

Choose randomly $m$ members of $A''$.  We shall prove that the probability 
that the set they form has exactly $m$ members and satisfies clause ($f$), 
is positive.  This suffices.  Let us explain.  We fix $n^*$ among these
elements and call the set they form $u$. 
\newline 

In clause $(*)_{12}$ for $\ell = 1,2$ we let $v_\ell =: \{i \in w_\ell:|u| 
+ |i \cap w_\ell| \in g_k(u)\}$. \newline
By $(*)_{15}$ the problem is that $\Rightarrow$ may fail. \newline

Let $x_1,\dotsc,x_m$ be random variables on $A''$.  The probability that
$\dsize \bigvee_{i \ne j}x_i = x_j$ is $\le \binom m2 \cdot 
{\frac 1{|A''|}}$. \newline

Now for $k < k(*)^1,u \in [\{1,\dotsc,m\}]^{n^*},w_1,w_2 \in [\{1,\dotsc,m\}
\backslash u]^{n^\otimes-n^*}$,\newline
$v_1 \subseteq w_1,v_2 \subseteq w_2$ defined as above, and a
possible linear order $<^*$ on $u^* = u \cup w_1 \cup w_2$, we shall give 
an upper bound for the probability that

$$
\dsize \bigwedge_{\ell_1,\ell_2 \in u^*} (\ell_1 <^* \ell_2 \Leftrightarrow
x_{\ell_1} < x_{\ell_2})
$$

\noindent
and they form a counterexample to clause ($f$) (of Claim 1.5).  So
in particular $u < w_1,u < w_2$. \newline
Choose $\ell \in v_1 \backslash v_2$ (as $v_1 \ne v_2$ and
$|v_1| = |g_k(u)| = |v_2|$ it exists).  We can
first draw $x_j$ for $j \ne \ell$.  Now we know $f_k(u \cup w_2)$;
note: we may not
know $w_2$ as possibly $\ell \in w_2$, but as $\ell \notin v_2$, 
by the choice of $A'$ it is not necessary to know $w_2$.  Now there is at
most one bad choice of $x_\ell$ (the others are good (inequality) or
irrelevant ($<^*$ is not right) by $(d)$ + $(e)$) so the probability of 
this is $\le {\frac 1{|A'|}}$.  So if we fix the set
$u^* = u \cup w_1 \cup w_2$ and concentrate on the case $|u^*| = 2n^\otimes
-n^*$, we have $2n^\otimes-n^*$ possibilities to choose $\ell \in u^*$ and
then having to choose $x_i$ for $i \ne \ell$, we know $u$ and have
$\le {\binom {2(n^\otimes-n^*)-1}{n^\otimes - n^*}}$ ways to choose $w_2$, 
so the probability of failure is
$\le (2n^\otimes - n){\binom {2(n^\otimes-n^*)-1}{n^\otimes - n^*}} \cdot
{\frac 1{|A''|}}$. \newline
So the probability that some failure occurs is at most (the cases
$|u^*| < 2n^\otimes - n^*$ and $x_1 = x_2$ are swallowed when $m$ is not too
small)

$$
(m^{2n^\otimes - n^*} \cdot k(*)(2n^\otimes - n)
{\binom {2(n^\otimes-n^*)-1}{n^\otimes - n^*}} \cdot
{\frac 1{|A''|}}.
$$
\medskip

Now by assumption $(*)_{13}$ this probability is $< 1$ so the conclusion
is clear. \hfill$\square_{1.5}$
\bigskip

Before we state and prove the last fact, which finishes the proof of the
theorem, we remind the reader of the following observation.  The proof is
easily obtained by induction on $\ell$.
\bigskip

\proclaim{Observation 1.6}  1) $\beth_\ell(kx) \ge k \beth_\ell(x)$ when
$x,k \ge 2$ and $\ell \ge 1$. \newline
2)  $\beth_\ell(kx) \ge (\beth_\ell(x))^k$ when $x \ge 2,k \ge 2$ and
$\ell \ge 1$.
\endproclaim
\bigskip

\proclaim{Fact 1.7}  Assume that $n^\otimes,n^*,m(n^*),k(*),\varepsilon,t$ and
$c_1$ are as in $(*)_9(a) \text{ and } (c)$. \endproclaim

Let us define

$$
c_2 = \text{ Max}\{(1 + \varepsilon)c_1,2\}
$$

$$
c_3 = n^\otimes \times (c_2)^2
$$
\medskip

\noindent
and the function $m(-)$ as follows: for $n \in (n^*,n^\otimes]$ by

$$
m(n) = \beth_{n-n^*}(m^{n^* + 1}c_3^{n-n^*})
$$
\medskip

\noindent
where

$$
m = m(n^*).
$$
\medskip

\noindent
\underbar{Then} $(*)_9(b)$ holds.
\bigskip

\demo{Proof}  We need to check that for $n \in [n^*,n^\otimes)$ 

$$
m(n+1) \ge 2^{(1 + \varepsilon)c_1m(n)^{n+1}}
$$
\medskip

\noindent
or equivalently

$$
\text{log}_2(m(n+1)) \ge (1 + \varepsilon)c_1 \, m(n)^{n+1}
$$
\medskip

\noindent
so it is enough that

$$
\text{log}_2(m(n+1)) \ge c_2 \, m(n)^{n+1}
$$
\medskip

\noindent
i.e.

$$
\text{log}_2(\beth_{n+1-n^*}(m^{n^* + 1} c^{n+1-n^*}_3)) \ge c_2 \, 
m(n)^{n+1}
$$
\medskip

\noindent
i.e., when $n > n^*$

$$
\beth_{n-n^*}(c_3(c^{n-n^*}_3) m^{n^*+1}) \ge
c_2(\beth_{n-n^*}(c^{n-n^*}_3 m^{n^*+1}))^{n+1}.
$$
\medskip

\noindent
It suffices by the above observation that

$$
\beth_{n-n^*}(c_3 \cdot c^{n-n^*}_3 m^{n^*+1}) \ge
\beth_{n-n^*}(c_2(n+1) c^{n-n^*}_3 m^{n^*+1}),
$$
\medskip

\noindent
which is true by the definition of $c_3$ when $n > n^*$.
\medskip

\noindent
For $n=n^*$ we need that

$$
m(n^*+1) \ge 2^{(1 + \varepsilon)c_1m^{n^*+1}}
$$
\medskip

\noindent
i.e.

$$
2^{m^{n^*+1}c_3} \ge 2^{(1 + \varepsilon)c_1m^{n^*+1}},
$$
\medskip

\noindent
which is true as $c_3 \ge c_2 \ge (1 + \varepsilon)c_1$.
\hfill$\square_{1.7}$
\enddemo
\bigskip

\proclaim{1.8 Proof of Lemma 1.1}  So $m,n,\varepsilon$ are given.  Let
\medskip
\roster
\item "{$(a)$}"  $n^*=1,n^\otimes = n,k(*) = 1,t=1,c_1$ as in 1.2, and
$c_2,c_3$ as in 1.7
\item "{$(b)$}"  $m_0 = m$
$$
\align
m_1 = &k(*) \cdot (2n^\otimes - n^*) \cdot {\binom{2(n^\otimes - n^*)-1}
{n^\otimes - n^*}}(m_0)^{2n^\otimes - n^*} \\
  =&(2n-1){\binom{2n-3}{n-1}} m^{2n-1}
\endalign
$$

$$
m_2 = (k(*)m_1)^{k(*)+1} = (m_1)^2
$$

$$
m_3 = (m_1)^2 \cdot 2^{k(*)} = 2(m_1)^2
$$
\item "{$(c)$}"  we define function $m(-)$ with domain $[n^*,n^\otimes]$:
\endroster

$$
\text{for } \ell = 1 \text{ we let } m(1) = m_3
$$

$$
\text{for } \ell > 1 \text{ we let } m(\ell) = 
\beth_{\ell-1}(c^{\ell-1}_3 \cdot (m_3)^2).
$$
\medskip

So we are given $m^\otimes > m(n)$.  In Claim 1.3 from the assumption
$(*)_9$, clauses $(a),(c)$ hold and clause $(b)$ holds by Fact 1.7.  Also
assumption $(*)_{10}$ of 1.3 holds (with $\{ 1,\dotsc,m^\otimes\}$ standing
for $A,f_0$ the given function $f$ (in 1.1), and $g$ constantly zero).

So there are $A \in [\{1,\dotsc,m^\otimes\}]^{m(1)+1},g',g_0,h_0$
satisfying the conclusion of 1.3 i.e. $(*)_{11}$.  $A$ here stands for
$A'$ in 1.3.  Note
$|A| = m_3 + 1$.  Now apply 1.4(2) with $A,m_3,m_2,f_0,g_0,1$ here standing
for $A,m(1),m(0),f_0,g_0,k(*)$ there and get $A' \in [A]^{m_2+1}$.  Next
we apply Claim 1.4(1) with $A',m_2,m_1,h_0,1$ here standing for $A',m(1),m(0),
h_0,k(*)$ there and get $A'' \in [A']^{m_1}$ satisfying the conclusion of
1.4(1) i.e. $(*)_{12}$.  Lastly apply Claim 1.5 and get $A^* \in [A'']^{m_0}
= [A'']^m$ satisfying the conclusion of 1.5; i.e. $(*)_{14}$.  Now $A^*$ is
as required.
\endproclaim
\bigskip

\remark{1.9 Remark}  We could have applied 1.5 in each stage, or just for
$n = 3$, this saves, somewhat.
\endremark

\newpage

REFERENCES
\bigskip

\bibliographystyle{literal-plain}
\bibliography{lista,listb,listx}

\def\germ{\frak} \def\scr{\cal}
  \ifx\documentclass\undefinedcs\def\rm{\fam0\tenrm}\fi
  \def\defaultdefine#1#2{\expandafter\ifx\csname#1\endcsname\relax
  \expandafter\def\csname#1\endcsname{#2}\fi} \defaultdefine{Bbb}{\bf}
  \defaultdefine{frak}{\bf} \defaultdefine{mathbb}{\bf}
  \defaultdefine{beth}{BETH} \def\bbfI{{\Bbb I}} \def\mbox{\hbox}
  \def\text{\hbox} \def\om{\omega} \def\Cal#1{{\bf #1}} \def\pcf{pcf}
  \defaultdefine{cf}{cf} \defaultdefine{reals}{{\Bbb R}}
  \defaultdefine{real}{{\Bbb R}} \def\restriction{{|}} \def\club{CLUB}
  \def\w{\omega} \def\exist{\exists} \def\se{{\germ se}} \def\bb{{\bf b}}
  \def\equivalence{\equiv} \let\lt< \let\gt> \def\cite#1{[#1]}
\begin{thebibliography}{GrRoSp}
\makeatletter \renewcommand{\@biblabel}[1]{[#1]} \makeatother

\bibitem[ErRa]{ErRa}Paul Erd\H{o}s and Richard Rado.
\newblock {A combinatorial theorem}.
\newblock {\em Journal of the London Mathematical Society}, {\bf 25}:249--255,
  1950.

\bibitem[ErSp]{ErSp}Paul Erd\H{o}s and Joel Spencer.
\newblock {\em {Probabilistic Methods in Combinatorics}}.
\newblock Academic Press, New York, 1974.

\bibitem[GrRoSp]{GrRoSp}Ronald Graham, Bruce~L. Rothschild, and Joel Spencer.
\newblock {\em {Ramsey Theory}}.
\newblock {Willey -- Interscience Series in Discrete Mathematics}. Willey, New
  York, 1980.

\bibitem[LeRo94]{LeRo94}Hanno Lefmann and Vojtech Rodl.
\newblock {???}

\bibitem[LeRo93]{LeRo93}Hanno Lefmann and Vojtech Rodl.
\newblock {On Canonical Ramsey Numbers for Complete Graphs versus Paths}.
\newblock {\em Journal of Combinatorial Theory, ser B}, {\bf 58}:1--13, 1993.

\bibitem[Ra86]{Ra86}Richard Rado.
\newblock {Note on canonical partitions}.
\newblock {\em Bulletin London Mathematical Society}, {\bf 18}:123--126, 1986.

\bibitem[Sh 37]{Sh:37}Saharon Shelah.
\newblock {A two-cardinal theorem}.
\newblock {\em {Proceedings of the American Mathematical Society}}, {\bf
  48}:207--213, 1975.

\end{thebibliography}

\shlhetal

\enddocument
\bye